\documentclass[10pt]{article}
\usepackage{latexsym,amsmath,graphics,graphicx}
\usepackage[psamsfonts]{amssymb}
\usepackage[bookmarks,bookmarksnumbered,pdfstartview=FitBH]{hyperref}

\newtheorem{theorem}{Theorem}[section]

\newtheorem{lemma}[theorem]{Lemma}
\newtheorem{proposition}[theorem]{Proposition}
\newtheorem{corollary}[theorem]{Corollary}
\newtheorem{definition}[theorem]{Definition}
\newtheorem{example}[theorem]{Example}
\newtheorem{remark}[theorem]{Remark}
\newtheorem{question}[theorem]{Question}

\def\endofproof{\nobreak\hfill $\blacksquare$ \goodbreak}

\newcounter{MEprop}
\newcounter{MEpropcourant}

\renewcommand{\theMEpropcourant}{\textbf{P}$_{\!\mbox{\tiny{ME}}} {\mathbf{\arabic{MEpropcourant}}}$}

\newenvironment{propME}
{\begin{list}
{\theMEpropcourant} 
{\usecounter{MEpropcourant}\leftmargin 50pt\rightmargin 40pt\labelsep 10pt}
\setcounter{MEpropcourant}{\value{MEprop}}
}
{\setcounter{MEprop}{\value{MEpropcourant}
}
\end{list}
\vskip-5pt}

\def\RR{\mathcal{R}}
\def\SS{\mathcal{S}}
\def\FF{\mathbf{F}}
\def\Nmath{\mathbb{N}}
\def\Zmath{\mathbb{Z}}
\def\Rmath{\mathbb{R}}
\def\Qmath{\mathbb{Q}}

\def\GG{\Gamma}
\def\LL{\Lambda}
\def\CC{C} 
\def\hhh{g}

\newcommand{\ME}{\raisebox{-0.5mm}{\ ${{\overset{\mathrm{ME}}{\sim}}}$\ }}
\newcommand{\notME}{\ {{\overset{\mathrm{ME}}{\not\sim}}}\ }
\newcommand{\MEa}[1][{}]{\ \underset{\smash{#1}}{{\overset{\mathrm{ME}}{\sim}}}\ } 
\def\OE{\raisebox{-0.5mm}{\ ${{\overset{\mathrm{OE}}{\sim}}}$\ }}

\newcommand{\OEstrong}{\underset{\smash{st}}{{\overset{\mathrm{OE}}{\sim}}}\ }
\newcommand{\itFP}[2]{\overset{#2}{\underset{\ #1}{\ast}}\ }

\newcommand{\amalgam}[1]{\underset{#1}{\ast}}

\def\cost{\mathcal{C}}
\def\pp{p}
\def\GGG{\mathcal{G}}
\def\AAA{A}
\def\BBB{B}
\def\bbb{b}
\def\aaa{a}
\def\XXX{X}
\def\YYY{Y}

\begin{document}

\thispagestyle{empty}

\title{Examples of Groups that are Measure Equivalent to the Free Group}
\author{D. Gaboriau\thanks{C.N.R.S.}}
\date{}
\maketitle
\begin{abstract}
{Measure Equivalence (ME) is the measure theoretic counterpart of quasi-isometry. This field grew considerably during the last years, developing tools to distinguish between different ME classes of countable groups. On the other hand, contructions of ME equivalent groups are very rare. We present a new method, based on a notion of measurable free-factor, and we apply it to exhibit a new family of groups that are measure equivalent to the free group. We also present a quite extensive survey on results about Measure Equivalence for countable groups.}
\end{abstract}



\section{Introduction}

During the last years, much progress has been made in the classification 
of countable discrete groups up to Measure Equivalence (ME). All of them 
describe criteria ensuring that certain groups don't belong to the same 
equivalence class (see for instance the work of R.~Zimmer \cite{Zim84,Zim91}, S.~Adams \cite{Ada90,Ada94,Ada95}, Adams-Spatzier \cite{AS90}, 
A.~Furman \cite{Fur99a,Fur99b}, D.~Gaboriau \cite{Gab00,Gab02}, 
Monod-Shalom \cite{MS02},...
We also would like to draw the attention of the reader 
to the very nice paper \cite{HK05} of 
 Hjorth-Kechris where similar results are developed.
 
 This notion (ME), introduced by M.~Gromov \cite{Gro93} is a measure 
theoretic analogue of Quasi-Isometry (QI) (see sect.~\ref{sect:generalities}). But, it takes its roots in the pioneering work of H.~Dye~\cite{Dye59,Dye63}, and 
even, somehow, back to Murray-von Neumann
\cite{MvN36}.

The simplest instances of ME groups are
 \textit{commensurable groups} or more generally \textit{commensurable up to finite kernel}\footnote{in the sense of \cite[IV.B.27]{Har00}: 
\textbf{commensurability up to finite kernel} is the equivalence relation on 
the set of countable groups generated by the equivalence of two terms of an 
exact sequence $1\to A\to B \to C\to 1$ as soon as the third term is a finite 
group. \label{foot: comm up to finite ker}}, and
 \textit{groups that are lattices} (=discrete, finite covolume subgroup)
 \textit{in the same locally compact second countable group}. 
 Recall that cocompact lattices are QI.

The first non-elementary iso-ME-class result is due to Dye \cite{Dye59,Dye63}
who puts together many amenable groups in the same ME-class, for example,
all the infinite groups with polynomial growth. A series of improvements 
led to Ornstein-Weiss' theorem \cite{OW80} asserting that, in fact, 
 \textit{all 
infinite amenable groups are ME to each other}. 
On the other hand, since amenability is 
a ME invariant, the ME class of $\Zmath$ consists precisely in 
all the infinite amenable groups \cite{Fur99a}.

Next, to build further ME groups, one can elaborate on these constructions via elementary procedures: direct products and free products 
(see Section~\ref{sect:generalities} for precise statements). 
And these are essentially the only known methods.

\bigskip
The family of groups ME to a free group\footnote{For consistency, we are led
to include the free groups $\FF_{\infty}$ on countably many elements and $\FF_{0}$ on $0$ element (i.e. the trivial group $\{1\}$).} contains for instance
the finite and amenable groups, as well as  the fundamental groups $\pi_1(\Sigma_g)$
of the compact orientable surfaces of genus $\geq 2$ (lattices in $\mathrm{SL}(2,\Rmath)$);
and it is stable under taking free products
(see Section~\ref{sect: free prod}, Property~\ref{P: free prod gps ME free gp})
and subgroups (see Section~\ref{sect: free prod}, Property~\ref{P: ME F stable subgroup}).
For instance, the following groups are ME to $\FF_2$:
$\FF_2 \ast \pi_1(\Sigma_g),\ \Zmath^2 \ast \Zmath^2,\ \Qmath\ast \Zmath/3\Zmath$, the triangular groups $T_{a,b,c}$ of isometries of the hyperbolic 
plane (with fundamental 
domain a triangle of angles $\frac{\pi}{a}, \frac{\pi}{b}, \frac{\pi}{c}$, 
$2\leq a\leq b\leq c$ and $\frac{1}{a}+\frac{1}{b}+\frac{1}{c}<1$),
$(\mathrm{SL}(2,\Zmath)\times \Zmath/5\Zmath)\ast \Qmath^2 \ast A$,
where $A$ is any amenable group.

Following a theorem of G.~Hjorth (see \cite{Hjo02} or also \cite[Th. 28.2, p. 98]{KM04}), being ME to a free group is
 equivalent with the following elementarily equivalent conditions
 (sect. \ref{sect: free prod}, \ref{P: Hjorth result}):\\
-- being treeable in the sense of Peres-Pemantle \cite{PP00}\\
-- admitting a treeable p.m.p.\footnote{p.m.p.: probability measure preserving, on a standard Borel probability space.} free\footnote{In this measure 
theoretic context, free means ``essentially free'', i.e. up to removing a set of measure zero.} action
in the sense of \cite{Ada90,Gab00}\\
-- having ergodic dimension 0 (for finite groups) or 1 in the sense of \cite{Gab02}

On the other hand, the free groups (and thus the family of groups ME to a free group)
split into four different ME-classes:
the class of $\FF_0$ (finite groups), that of $\FF_1=\Zmath$ (amenable groups),
that of $\FF_{p}$, $2 \leq p <\infty$, (all the $\FF_p$ are commensurable) and that
of $\FF_{\infty}$. The last two classes are distinguished by their $\ell^2$-Betti numbers -- $\beta_1\in (0,\infty)$, resp. $\beta_1=\infty$ -- \cite{Gab02} 
(or also by their cost -- $\cost \in (1,\infty)$, resp. $\cost=\infty$ -- \cite{Gab00}) (see Section~\ref{sect:generalities}).

However, the classification of groups ME to a free group seems, nowadays, 
completely out of reach\footnote{To compare with, recall that a group 
quasi-isometric to a free group is virtually a free group \cite{Sta68}.}.
In fact, the only groups whose ME-class is classified 
are finite groups, amenable groups and lattices 
in simple connected Lie groups with finite center and real rank $\geq 2$ \cite{Fur99a} (see Property \ref{P: ME rig. latt. high rank}).

\bigskip
\noindent
The contribution of this paper consists in a new construction of ME groups, leading to the exhibition of new groups that are ME to free groups $\FF_2$ and $\FF_{\infty}$.
\bigskip

Let $\GG$ be a countable group and $\LL$ a subgroup.
We denote by $\itFP{\LL}{n}\GG$
the \textbf{iterated amalgamated free product}
$\GG\amalgam{\LL}\GG\amalgam{\LL}\cdots \amalgam{\LL}\GG$ 
of $n$ copies of $\GG$ above the corresponding $n$
copies of $\LL$, where the injection morphisms are just the identity.

Consider the free group $\FF_{2p}=\langle a_{1},\cdots,a_{p},b_{1},\cdots,b_{p}\rangle$
and its cyclic subgroup $\CC$
generated by the product of commutators:
$\kappa:=\prod 
\limits_{i=1}^{i=p}[a_{i},b_{i}]$.
\begin{theorem}[Cor.~\ref{cor: branched surfaces ME F2}]
\label{th: branched surfaces ME F2}
For each $n\in \Nmath$, the iterated amalgamated free product
$\itFP{\CC}{n}\FF_{2p}$
is measure equivalent to the free group $\FF_2$.
\end{theorem}
This iterated amalgamated free product
$\itFP{\CC}{n}\FF_{2p}$ is the fundamental group
$\pi_1(\Sigma)$ of a ``branched surface'' (figure~1).
\begin{figure}[htbp]
\centering
\includegraphics[width=.30\textwidth]{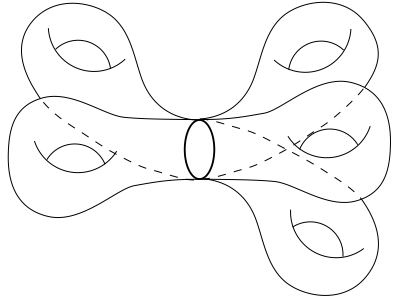}
\caption{Branched surface}
\end{figure}
Our examples for various $n$ are not obviously lattices in the same l.c.s.c. 
group (i.e. the question asked to several specialists remains open).
However, they are certainly not cocompact lattices in the same l.c.s.c. group nor commensurable up to finite kernel since they are not quasi-isometric:
\begin{proposition}
The number $n$ of factors $\FF_{2p}$ in the 
iterated amalgamated free product
$\itFP{\CC}{n}\FF_{2p}$ is a quasi-isometry invariant.
\end{proposition}
\textsc{Proof}: The boundary of the group is disconnected by pairs of points 
into either 
$1$, $2$ or $n$ connected components.\endofproof

\medskip
Theorem~\ref{th: branched surfaces ME F2} extends to an infinite amalgamated 
free product: $\itFP{\CC}{\infty}\FF_{2p}$ is measure equivalent to the free 
group $\FF_{\infty}$ (Cor.~\ref{cor: branched surfaces ME F2}).

\begin{corollary}
If $r$ is any non trivial element of $\Zmath$, then the amalgamated free product 
${\Gamma_r=\Zmath\amalgam{r=\kappa}\FF_{2p}}$ is measure equivalent to the free 
group $\FF_2$.
\end{corollary}
\textsc{Proof}: Consider the natural homomorphism $\smash{\Zmath\amalgam{r=\kappa}\FF_{2p}\to \Zmath/r\Zmath}$ that kills $\FF_{2p}$. Its kernel has index $r$ in $\Gamma_r$ and is isomorphic with $\itFP{\CC}{r}\FF_{2p}$.\endofproof
More generally, one gets:
\begin{theorem}(Cor.~\ref{cor: f.p. amalg. ME free group})
Let $G$ be any countable group, $H$  an infinite cyclic subgroup
and $\CC$ the cyclic subgroup generated by $\kappa$ in $\FF_{2p}$. Then 
$\smash{G\amalgam{H=\CC}\FF_{2p}}$ is ME to $G\amalgam{}\FF_{2p-1}$.
In particular, if $G$ is ME to $\FF_2$, then $G\amalgam{H=\CC}\FF_{2p}$ is ME to $\FF_2$.
\end{theorem}
By contrast, an amalgamated free product $\Zmath\amalgam{\Zmath}\Zmath$
is measure equivalent to a free group
if and only if 
the common $\Zmath$ coincides with one the 
components or injects as a subgroup of index 2 in both
(see \ref{P: l2 inv ME}**).
Also, if instead of $\CC$ one considers $2\CC$, generated by $\kappa^{2}$, 
and if $H$ is contained as a subgroup of index at least three in a greater 
abelian group $H'$ of $G$, then $G\amalgam{H=2\CC}\FF_{2p}$
cannot be ME to a free group (by \ref{P: ME F stable subgroup} and \ref{P: l2 inv ME}**):
it contains the nonamenable subgroup $H'\amalgam{H=2\CC}\CC$
with $\beta_1=0$.

\bigskip
To obtain these results we show that the subgroup $C$ happens to appear
as a free factor of the group $\FF_{2p}$ in a measure theoretic sense 
(see Def.~\ref{def: measure free-factor}). 
Our proof  makes use of standard percolation techniques
developed in \cite{LP05}.

\bigskip
\noindent
\textbf{Acknowledgment}:
I would like to thank Russell Lyons for worthwhile discussions.
I am deeply grateful to Yuval Peres for crucial indications in proving
Theorem~\ref{th: prod commutators is a measure-free-factor}.
Contrarily to his opinion, his contribution is really substantial.
I'm also grateful to Roman Sauer for drawing my attention to the co-induced action construction, and together with Yves de Cornulier, Nicolas Monod and Yehuda Shalom for useful comments on a preliminary version.

\section{Generalities about Measure Equivalence}
\label{sect:generalities}

\begin{definition}[{\cite[0.5.E]{Gro93}}]
Two countable groups $\Gamma_1$ and $\Gamma_2$ are \textbf{Measure
Equivalent (ME)}
if there exist commuting, measure preserving, free actions of
$\Gamma_1$ and $\Gamma_2$ on some Lebesgue measure space
$(\Omega, m)$ such that the action of each of the groups
admits a {\em finite measure} fundamental domain
($D_i$ for $\Gamma_i$). 
\end{definition}
In this case we say that $\Gamma_1$ is ME to $\Gamma_2$ with \textbf{index}
$\iota=[\Gamma_1:\Gamma_2]_{\Omega}:=m(D_2)/m(D_1)$,
and we denote the existence of such a \textbf{coupling} $\Omega$
by
$$\Gamma_1\MEa[\iota] \Gamma_2\mathrm{\ \ \ or\ simply\ \ \ } \Gamma_1\ME\Gamma_2$$ 
when the particular value of the index is irrelevant.
Measure Equivalence is an equivalence relation on the set of countable groups.
The grounds are established by A.~Furman in \cite[Sect. 2]{Fur99a}.

\medskip

The notion of Measure Equivalence 
has been introduced by M.~Gromov as the 
\textit{measure theoretic} analogue of the 
notion of Quasi-Isometry, which is more 
\textit{topological}, as testified by the following criterion:

\medskip
\noindent
\textbf{Criterion for quasi-isometry (\cite[0.2.C$'_{2}$]{Gro93})}
Two finitely generated groups $\Gamma_1$ and $\Gamma_2$ are
{\em quasi-isometric (QI)}\index{quasi!isometry}\index{isometry!quasi---}\index{QI} 
i{f}{f} there exist commuting,
continuous actions of
$\Gamma_1$ and $\Gamma_2$ on some locally compact space $M$,
such that the action of each of the groups is properly discontinuous
and has a {\em compact} fundamental domain.

\medskip
We anthologize basic and less basic properties of ME, with references, or proofs when necessary.

\subsection{Basic Properties of ME}

\begin{propME}
\item The ME class of the trivial group $\{1\}$ consists in 
	all finite groups. \label{P: trivial group}
\item Groups that are commensurable up to finite kernel are ME.
\label{P: commensur}
\end{propME}
Commensurability  up to finite kernel 
(footnote~{\footnotesize\ref{foot: comm up to finite ker}})
 corresponds exactly to ME
with a countable atomic space~$\Omega$.
\begin{propME}
\item Lattices in the same locally compact second countable group are ME.
\end{propME}
Actions of the two lattices by left multiplication and right multiplication by the inverse on the group deliver the desired ME coupling.
\begin{propME}
\item (Direct Product)
	$\Gamma_1 \ME \Lambda_1$ and $\Gamma_2\ME \Lambda_2$ $\Rightarrow$
	$\Gamma_1\times \Gamma_2\ME \Lambda_1\times \Lambda_2$
\end{propME}
This is quite obvious, by taking the product actions on the product space.
The converse fails in general. For instance by Ornstein-Weiss' theorem
recalled in the introduction (see also \ref{P: amen ME inv}), all infinite amenable groups are ME so that\footnote{More generally $(\Gamma\times A)\ME(\Gamma\times A)\times B$, for any group $\Gamma$, and any infinite amenable groups $A$ and $B$.
} $\FF_2\times \Zmath\ME \FF_2\times (\Zmath\times\Zmath)=(\FF_2\times \Zmath)\times\Zmath$.
but 
$\FF_2\notME\FF_2\times \Zmath$ (see \cite[Th. 6.1]{Ada94} or
\ref{P: l2 inv ME}**).

However, it becomes true when restricted to a quite large class of groups
(see \ref{P: prod gps in C}). 

\medskip
Recall that two p.m.p. actions $\alpha_i$ of $\Gamma_i$ on $(X_i,\mu_i)$ ($i=1,2$)
are \textbf{stably orbit equivalent (SOE)}
if there are complete sections $A_i\subset X_i$
such that the orbit equivalence relations
$(\RR_{\alpha_i(\Gamma_i)}, \mu_i)$ restricted to $A_i$ are isomorphic 
via a measure scaling  
isomorphism $f:A_1\to A_2$ ($f_*\mu_{1\vert A_1}=\lambda\, \mu_{2\vert A_2}$).
The \textbf{compression constant} is defined as $\iota:=\mu_2(A_2)/\mu_1(A_1)=1/\lambda$.
The connection with ME is very useful 
and is proved in \cite[Lem. 3.2, Th. 3.3]{Fur99b}, where it is credited to
Zimmer and Gromov.
An easy additional argument gives the freeness condition
(see \cite[Th. 2.3]{Gab00b}).
\textbf{Orbit Equivalence (OE)} corresponds to the case where $A_i=X_i$ 
(and thus $\iota=1$)
and where the fundamental domains $D_1,D_2$ of the measure equivalence coupling may be taken
to coincide.  

\begin{propME}
\item (ME $\Leftrightarrow$ SOE) \label{P: ME <-> SOE}
-- Two groups $\Gamma_1$ and $\Gamma_2$ are ME (with index $\iota$) 
if and only if they admit
stably orbit equivalent \textbf{(SOE)} free p.m.p. actions 
(with compression constant~$\iota$).\\
-- They are ME with a common fundamental domain if and only if 
they admit orbit equivalent \textbf{(OE)} free p.m.p. actions.
We shall denote this last situation by $$\Gamma_1\OE\Gamma_2$$
-- If the only possible value of the index $\iota=[\Gamma_1:\Gamma_2]$
is $1$, then $\Gamma_1\OE\Gamma_2$.
\end{propME}
For its conciseness, we present a proof of this result.\\
\textsc{Proof:} ME $\Rightarrow$ SOE.
Consider an action of $\Gamma_1\times \Gamma_2$ on $({\Omega},m)$ witnessing $\Gamma_1\MEa[\iota]\Gamma_2$.
The existence of fundamental domains of finite measures $D_i$ ensures the existence of 
a $\Gamma_1$-equivariant Borel map onto a finite measure Borel space 
$\pi_1: \Gamma_1\times \Gamma_2  \curvearrowright ({\Omega},m) \to 
\Gamma_1 \curvearrowright X_1:={\Omega} \slash \Gamma_2\simeq D_2$
and similarly for $\pi_2$ by exchanging the subscripts $1$ and $2$.
Any two points $y,y'$ of $\Omega$ are $\Gamma_1\times \Gamma_2$-equivalent
i{f}{f} 
their images under $\pi_i$ are $\Gamma_i$-equivalent ($i=1,2$).
Since $\pi_1,\pi_2$ have countable fibers, one can find a non-null ${Y}\subset \Omega$ on which 
both $\pi_1$ and $\pi_2$ are one-to-one. The natural isomorphism $f$ between $\pi_1(Y)$ and $\pi_2(Y)$
preserves the restriction of the equivalence relations. Normalizing the induced measure on $X_i$
to a probability measure leads to identify index and 
compression constant. The freeness is obtained by just considering any free p.m.p. action of $\Gamma_1$
on some space $(Z,\nu)$ (with trivial action of $\Gamma_2$)
and replacing the action on $\Omega$ by the diagonal action on the product 
measure space $\Omega\times Z$. In case of coincidence, $Y:=D_1=D_2\simeq X_1\simeq X_2$ 
delivers OE actions.

\noindent
\textsc{Proof:} SOE $\Rightarrow$ ME.
Let $\Gamma_i$ act on $(X_i,\mu_i)$ and assume $X_1\supset {Y}_1\simeq {Y}_2 \subset X_2$ witnesses 
the SOE. 
Rescale the measures $\mu_1, \mu_2$ so that ${Y}_1$ and ${Y}_2$
get the same measure and so that the glueing $X:= (X_1\coprod X_2)\slash ({Y}_1={Y}_2)$
 of $X_1,X_2$ along $Y_1,Y_2$ becomes a probability space.
Define  on $X$ the equivalence relation $\RR$ extending $\RR_{\Gamma_i}$ on $X_i$.
Consider the coordinate projections $\pi_1$ and $\pi_2$ from $\RR$ to~$X$.
The Borel subset $\Omega\subset \RR\subset X\times X$ defined as the preimage by $\pi_1\times \pi_2$ of
$X_1\times X_2$ gives the desired ME coupling:
$\Gamma_i$  acts on the $i$-th coordinate (and trivially on the other one).
Any Borel section of $\pi_1$ gives a Borel fundamental domain for $\Gamma_2$
and any Borel section of $\pi_2$ gives a Borel fundamental domain for $\Gamma_1$.
Their measures, for the natural measure on $\RR$ (see \cite{FM77a}), are compared via the rescaling
to the compression constant. In case of OE, $\Omega=\RR$ and the diagonal $\Delta\subset \RR$
is a fundamental domain for both groups.
Observe that ergodic SOE actions on $X_i$ are in fact OE if and only if
$\iota=1$. The condition, when 
$\iota=1$ is the only possible index, follows by considering an ergodic component. \endofproof

\subsection{Free Products and Amalgamations}\label{sect: free prod}
Stability under taking free products requires the additional assumption
that not only the groups are ME but also that they admit a coupling with a common 
fundamental domain.
That such a result may be true is suggested in \cite[Rem. 2.27]{MS02}.
We prove it below (after Lemma~\ref{lem: loc 1-to-1 and splitting}) as well as some generalizations.
\begin{propME}
\item (Free Product) $\Gamma_1\OE\Lambda_1$ and $\Gamma_2\OE\Lambda_2$\ $\Rightarrow$ \ 
	$\Gamma_1\ast \Gamma_2\OE\Lambda_1\ast \Lambda_2$\label{P: stability under free products}
\end{propME}
And more generally:
\begin{propME}
\item[\theMEpropcourant $^{*}$] (Infinite Free Product) 
	$\Gamma_j\OE\Lambda_j$ for all $j\in \Nmath$
	\ $\Rightarrow$ \ 
	$\displaystyle{\underset{j\in\Nmath}{\ast}}\Gamma_j\OE
	\displaystyle{\underset{j\in\Nmath}{\ast}}\Lambda_j$
\end{propME}
The hypothesis $\OE$ cannot be replaced by 
$\ME$ since the following three groups belong to three different ME-classes\footnote{\label{footn: Papasaglu-Whyte} The same phenomenon occurs in the analogous setting of geometric group theory and for similar reasons: Bi-Lipschitz equivalence of groups passes to free products, but not just QI of them. However, the (counter-)example of $\{1\}\ast\Zmath/2\Zmath$, $\Zmath/2\Zmath\ast\Zmath/2\Zmath$ and $\Zmath/3\Zmath\ast\Zmath/2\Zmath$ is essentially the only one and is satisfactorily worked around by a result of Papasoglu-Whyte: \textit{If two families of non trivial ($\not=\{1\}$) finitely generated groups $(\Gamma_i)_{i=1,\cdots, s}$ and $(\Lambda_j)_{j=1,\cdots, t}$ $s,t \geq 2$ define the same sets of quasi-isometry types (without multiplicity), then the free products $\Gamma=\Gamma_1\ast\Gamma_2\ast\cdots \Gamma_s$ and $\Lambda=\Lambda_1\ast \Lambda_2 \ast\cdots \ast\Lambda_t$ are quasi-isometric unless $\Gamma$ or $\Lambda$ is $\Zmath/2\Zmath \ast \Zmath/2\Zmath$ \cite[Th. 02]{PW02}.}
}:
$\{1\}\ast\Zmath/2\Zmath$ (finite) (see \ref{P: trivial group}), $\Zmath/2\Zmath\ast\Zmath/2\Zmath$ (amenable) (see \ref{P: amen ME inv}) and $\Zmath/3\Zmath\ast\Zmath/2\Zmath$ (nonamenable). Similarly, free groups of different ranks are not OE (for a proof, we draw the attention of the reader to the simplified preliminary version \cite{Gab98} of \cite{Gab00}), and the groups $\FF_p\ast (\FF_2\times \FF_2)$ all belong to different ME-classes for different $p$'s: their $\ell^{2}$ Betti numbers\footnote{The sequence of $\ell^2$ Betti numbers of $\FF_p\ast (\FF_2\times \FF_2)$ is $(\beta_i(\FF_p\ast (\FF_2\times\FF_2)))_{i\in\Nmath}=(0, p, 1, 0, 0, \cdots, 0, \cdots)$. Incidentally, no one of QI or ME classification is finer than the other one, as illustrated by this family of groups (that they are QI is due to K.~Whyte, see footnote~\ref{footn: Papasaglu-Whyte}).} do not agree with \ref{P: l2 inv ME}. However, this hypothesis may be relaxed when working with groups ME to free groups:
\begin{propME}
\item If each $\Gamma_j$ is ME to some free group, then $\underset{j\in\Nmath}{\ast}\Gamma_j$ is ME to a free group.\label{P: free prod gps ME free gp}
\end{propME}
In particular,
\begin{propME}
\item[\theMEpropcourant $^{*}$] $\Gamma_1\ME\FF_p$
	and $\Gamma_2\ME\FF_q$, for $p,q\in\{1,2,\infty\}$ \  $\Rightarrow$ \ $\Gamma_1\ast \Gamma_2
	\ME\FF_{p+q}$
\end{propME}
For instance, non trivial free products of finitely many infinite 
amenable groups are ME to $\FF_2$ (see \cite[p. 145]{Gab00}).
Observe that free products of finitely many finite groups are virtually free.
\footnote{It follows from \ref{P: free prod gps ME free gp} that all the groups $\Gamma_{S}$ that appear in \cite[Th.~1]{HK05}
are ME when $S$ is finite (resp. infinite).}

\bigskip

We make an intensive use of countable measured equivalence relations $\RR$
on standard probability measure spaces (see \cite{FM77a} and the reminders of 
\cite{Gab00, Gab02}) : \textbf{they all are assumed to preserve the measure}.

A countable group $\Gamma$ is \textbf{treeable} in the sense of Peres-Pemantle
\cite{PP00} if the set of trees with vertex set $\Gamma$ supports a
$\Gamma$-invariant probability measure.
The group
$\Gamma$ is \textbf{not anti-treeable} in the sense of \cite{Gab00}
if it admits a treeable free action $\alpha$, i.e. an action
whose associated equivalence relation $\RR_{\alpha}$ 
is generated by a treeing. Recall that a \textbf{graphing} of $\RR$
is a countable family $\Phi=(\varphi_j)_{j\in J}$ of partially defined
isomorphisms $\varphi_j:A_j\to B_j$ between Borel subsets of $X$
such that $\RR$ is the smallest equivalence relation 
satisfying for all $j\in J$ for all $x\in A_j$,  
$x\sim \varphi_j(x)$. It equips each orbit with a graph structure.
When these graphs are trees, the graphing is called a \textbf{treeing}
\cite{Ada90} and the equivalence relation \textbf{treeable}.

The two following notions taken from \cite{Gab02}, which we mention just for completeness (see also \ref{P: ergodic dimension ME inv.}),
 won't be seriously used in the remainder of the paper.
Equivalence relations, as groupoids, admit discrete actions on fields of simplicial complexes.
The \textbf{ergodic dimension} is the smallest possible dimension
of such a field of contractible simplicial complexes.
For the \textbf{approximate ergodic dimension}, one considers increasing exhaustions
of $\RR$ by measurable subrelations and the smallest possible ergodic dimensions
of these approximations.
\begin{propME}
\item \label{P: Hjorth result}
The following conditions on a countable group $\Gamma$ are equivalent:\\
\phantom{.}\hskip5pt 
-i- being treeable in the sense of Peres-Pemantle\\
\phantom{.}\hskip5pt
-ii- admitting a treeable p.m.p. free action
in the sense of \cite{Gab00}\\
\phantom{.}\hskip5pt 
-iii- having ergodic dimension 0 (for finite groups) or 1\\
\phantom{.}\hskip5pt
-iv- being ME to a free group
\end{propME}
The equivalence of the first three conditions is elementary, while the
connection with the fourth one is a result of G.~Hjorth:
see \cite{Hjo02} or also \cite[Th. 28.2, Th. 28.5]{KM04}. In order
to apply \cite[Th. 28.2, Th. 28.5]{KM04}, observe that
after considering an ergodic component, the cost can be arranged, by SOE, 
to be an integer.

\bigskip
We introduce some terminology.
Let $\RR$ and $\SS$ be two countable measured equivalence relations 
on the standard probability measure spaces $(X,\mu)$ and $(Y,\nu)$ respectively.
A measurable map $\pp:X\to Y$ is a \textbf{morphism}
from $\RR$ to $\SS$ if almost every $\RR$-equivalent points $x,x'$ of $X$
have $\SS$-equivalent images $\pp(x), \pp(x')$ and if the pushforward measure 
$p_*\mu$ is equivalent to $\nu$.
Say that such a morphism is \textbf{locally one-to-one}
if, for almost every $x\in X$, $\pp$ induces a one-to-one map from 
the $\RR$-class of $x$ to the $\SS$-class of $\pp(x)$.
\begin{example}\label{ex: p loc 1-to-1}
If $\RR_{\alpha}$ and $\SS_{\beta}$ are the orbit equivalence relations of two actions $\alpha,\beta$ of a 
countable group $\Gamma$ on $X$ (resp. $Y$), then any $\Gamma$-equivariant 
measurable map $\pp:X\to Y$ is a morphism from $\RR_{\alpha}$ to $\SS_{\beta}$. If the action $\beta$ is moreover free, then $\pp$ is locally one-to-one.
\end{example}
\begin{lemma}\label{lem: loc 1-to-1 -> action}
If $p:X\to Y$ is a locally one-to-one morphism from $\RR$ to $\SS$
and if $\SS$ is the orbit equivalence relation of an action $\beta$
of a countable group $\Gamma$, then there exists a unique action $\alpha$
of $\Gamma$ on $X$ inducing the equivalence relation $\RR$ and for which $\pp$
is $\Gamma$-equivariant. If $\beta$ is moreover free, then $\alpha$ is free.
\endofproof
\end{lemma}
And more generally, for graphings (see \cite{Gab00}): 
\begin{lemma}\label{lem: loc 1-to-1 -> graphing}
If $p:X\to Y$ is a locally one-to-one morphism from $\RR$ to $\SS$
and if $\SS$ is generated by a graphing $\Phi=(\varphi_{1}, \varphi_{2}, \cdots, \varphi_{j}, \cdots)$, then there exists a unique graphing 
$\tilde{\Phi}=(\tilde{\varphi}_{1}, \tilde{\varphi}_{2}, \cdots, 
\tilde{\varphi}_{j}, \cdots)$ generating $\RR$ and for which $p$ is equivariant (i.e. $\forall \varphi_j$,
$x\in \mathrm{Dom}(\tilde{\varphi}_j)$ i{f}{f} $p(x)\in \mathrm{Dom}(\varphi_j)$ and in this case
$p \tilde{\varphi}_{j}(x)=\varphi_j(p(x))$).
In particular, if $\SS$ is treeable, then $\RR$ is treeable.
\endofproof
\end{lemma}

We have introduced in \cite[Sect. IV.B]{Gab00} the notions of an equivalence relation $\RR$
that splits as a \textbf{free product} of two (or more) subrelations $\RR=\RR_1\ast \RR_2$ and more generally  of an equivalence relation $\RR$
that splits as a \textbf{free product with amalgamation} 
of two subrelations over a third one
$\RR=\RR_1\amalgam{\RR_3} \RR_2$:\\
-- $\RR$ is \textit{generated} by $\RR_1$ and $\RR_2$ (i.e. $\RR$ is the smallest equivalence relation containing $\RR_1$ and $\RR_2$ as subrelations)\\
-- for (almost) every $2p$-tuple $(x_j)_{j\in \Zmath/2p\Zmath}$ that are successively $\RR_1$ and $\RR_2$-equivalent (i.e. $x_{2i-1}\overset{\RR_{1}}{\sim} x_{2i} \overset{\RR_{2}}{\sim} x_{2i+1}$ $\forall i$), there exists two consecutive ones that are $\RR_3$-equivalent (i.e. 
there exists a $j$ such that $x_{j}\overset{\RR_{3}}{\sim} x_{j+1}$).

Free products correspond to triviality of $\RR_3$, i.e. the second condition becomes: there exists a $j$ such that $x_{j} = x_{j+1}$.

Standard examples are of course given by 
free actions of groups that split in a similar way.
\begin{lemma}\label{lem: split and group action}
Assume that $\RR$ splits as a free product $\RR=\RR_1\ast\RR_2$
and that $\RR_i$ is produced by a free action of $\Gamma_i$.
Then the induced action of the free product $\Gamma_1\ast \Gamma_2$ 
produces $\RR$ and is free.
\endofproof
\end{lemma}
\begin{lemma}\label{lem: loc 1-to-1 and splitting}
If $p:X\to Y$ is a locally one-to-one morphism from $\RR$ to $\SS$
and if $\SS$ splits as a free product with amalgamation
$\SS=\SS_1\amalgam{\SS_3} \SS_2$, then $\RR$ admits a corresponding splitting
$\RR=\RR_1\amalgam{\RR_3} \RR_2$ for which $p$ is 
a locally one-to-one morphism from $\RR_i$ to $\SS_i$, for $i=1,2,3$.
\endofproof
\end{lemma}

\bigskip
\noindent
\textsc{Proof} of \ref{P: stability under free products}$^{*}$
(the proof of \ref{P: stability under free products} follows by letting 
$\Gamma_j=\Lambda_j=\{1\}$ for $j\geq 3$):
Let $\SS_j$, ($j\in \Nmath$), be measured equivalence relations on 
$(Y_j,\nu_j)$, given by orbit equivalent free actions of $\Gamma_j$ and $\Lambda_j$.
Take any free p.m.p. action of 
${\ast}_{j\in\Nmath}\Gamma_j$ on a standard probability measure space $(Z,\nu)$
and consider the diagonal action on the product measure space 
 $(X,\mu):=(Z\times \prod_{j\in\Nmath} Y_j,\nu\times \prod_{j\in\Nmath}\nu_j)$
(i.e. $\gamma.x=\gamma.(z,y_1,y_2,\cdots,y_j,\cdots)=(\gamma.z,\gamma.y_1,\gamma.y_2,
\cdots,\gamma.y_j,\cdots)$, where ${\ast}_{j\in\Nmath}\Gamma_j$ acts
on $Y_i$ via the natural homomorphism ${\ast}_{j\in\Nmath}\Gamma_j\to \Gamma_i$).
This action is free and p.m.p. Denote by $\RR$ the equivalence relation it defines on $X$ and by $\RR_j$ the equivalence sub-relations defined by
the restriction of this action to $\Gamma_j$.
The measurable map $p_i:X=Z\times \prod_{j\in\Nmath} Y_j\to Y_i$
is $\Gamma_i$-equivariant and thus is a locally one-to-one morphism from 
$\RR_i$ to $\SS_i$ (ex.~\ref{ex: p loc 1-to-1}).
Lemma~\ref{lem: loc 1-to-1 -> action} provides
a free action of $\Lambda_i$ on $X$ producing $\RR_i$; and thus an 
action of the free product ${\ast}_{j\in\Nmath}\Lambda_j$ giving $\RR$.
The free product structure $\RR=\RR_1\ast \RR_2\ast\cdots\ast \RR_j\ast\cdots$ of the equivalence relation, given by the free product structure of the group ${\ast}_{j\in\Nmath}\Gamma_j$,
(see \cite[D\'ef. IV.9]{Gab00}) ensures that the action of the free product 
${\ast}_{j\in\Nmath}\Lambda_j$ is free (see Lemma~\ref{lem: split and group action}).
We have constructed orbit equivalent free actions of 
${\ast}_{j\in\Nmath}\Gamma_j$ and ${\ast}_{j\in\Nmath}\Lambda_j$.
\endofproof

\bigskip
\noindent
\textsc{Proof} of \ref{P: free prod gps ME free gp}:
While following the lines of the preceding proof, forget about the groups 
$\Lambda_j=\FF_{p_j}$ and just retain that $\SS_j$, generated by a free action of $\Gamma_j$, is treeable. The treeings of $\RR_j$ ($j\in \Nmath$), delivered
by lemma~\ref{lem: loc 1-to-1 -> graphing}, together define a graphing of 
$\RR=\RR_1\ast \RR_2\ast\cdots\ast \RR_j\ast\cdots$ that happen to be a treeing due to the free product structure of $\RR$. Then apply \ref{P: Hjorth result}.\endofproof

\begin{propME}
\item The family of groups ME to some free group is stable under taking subgroups.
\label{P: ME F stable subgroup}
\end{propME}
This follows from \cite[Th. 5]{Gab00} and \ref{P: Hjorth result}.

\subsection{Some Invariants}

\begin{propME}

\item (Amenability) The ME class of $\Zmath$ consists in 
	all infinite amenable groups.
Moreover, any ergodic action of any two infinite amenable groups are OE
and also SOE for any compression constant $\iota$ \cite{OW80}.
	\label{P: amen ME inv}
\end{propME}
This property is sometimes 
interpreted as a kind of \textit{elasticity}. This in wide contrast
with the \textit{rigidity} phenomenon described by Zimmer and Furman 
\cite{Zim84,Fur99a,Fur99b} for lattice in higher rank semi-simple Lie groups see also~\ref{P: ME rig. latt. high rank}.

Recall that Kazhdan property T is often considered as opposite to amenability.
More recently, the general opinion switched, after M.~Gromov to oppose Kazhdan 
property T with the greater class of groups with Haagerup property 
(= a-T-menability, see \cite{CCJJV01}).
\begin{propME}
\item (Kazhdan Property T) Kazhdan property T is a ME invariant \cite{Fur99a}.
	\label{P: T is ME inv}
\end{propME}
Lattices of $\mathrm{Sp}(n,1)$ are not ME to lattices of 
$\mathrm{SU}(p,1)$ or $\mathrm{SO}(p, 1)$.
We already knew from \cite{AS90} that groups 
with property T are not ME to groups that split non trivially as 
free product with amalgamation or HNN-extension.
\begin{propME}
\item (Cowling-Haagerup constant) Cowling-Haagerup constant $\Lambda(\Gamma)$ is a ME invariant \cite{CZ89,Jol01}.\label{P: C(Lambda) is ME inv}
\end{propME}
The Cowling-Haagerup constant $\Lambda(\Gamma)$ of the group $\Gamma$ 
is the infimum of the constants $C$ such that the Fourier algebra
$\mathcal{A}(\Gamma)$ admits an approximate unit bounded by $C$ in the
multipliers norm.
Jolissaint in \cite{Jol01} extends the above result of Cowling-Zimmer 
\cite{CZ89} from OE to ME.
Lattices in various $\mathrm{Sp}(n,1)$ are not ME ($\Lambda(\Gamma)=2n-1$).
Recall that $\Lambda(\Gamma)=1$ implies that $\Gamma$ has Haagerup property,
but that it is unknown whether the converse holds (Cowling conjecture).
\begin{propME}
\item (a-T-menability) Haagerup property is a ME invariant 
\end{propME}
This result was obtained independently by several people, including 
Jolissaint\cite{Jol01}, Popa \cite{Pop01} and Shalom (personal communication).
\medskip

Two kinds of numerical invariants of countable groups are
studied in \cite{Gab00,Gab02} in connection with ME:
Cost $\cost(\Gamma)\in [0,\infty]$ and $\ell^2$ Betti numbers 
$\beta_0(\Gamma), \beta_1(\Gamma), \cdots,  \beta_i(\Gamma), \cdots$.
\begin{propME}
\item (Cost) 
Having cost $=1$ (resp. $=\infty$) is a ME invariant \cite{Gab00}.
\label{P: cost}
\end{propME}
For all the groups for which the computation 
has been carried out until now, $\cost(\Gamma)=\beta_1(\Gamma)-\beta_0(\Gamma)+1$.
\begin{propME}
\item ($\ell^2$ Betti Numbers) ME groups $\Gamma$ and $\Lambda$ have
	proportional $\ell^2$ Betti numbers:
	there is $\lambda>0$ such that for every $n\in \Nmath$,
	$\beta_n(\Gamma)=\lambda \beta_n(\Lambda)$ \cite{Gab02}.
	\label{P: l2 inv ME}
\end{propME}
This gives\footnote{For $p\not=0$, $\beta_1(\FF_p)=p-1$ and all the other $\beta_i$ are zero} the splitting of groups ME to free groups into four ME-classes
according to whether $\beta_0\not=0, \beta_1=0$ (finite groups)
$\beta_0=0,\ \beta_1=0$ (amenable ones),
$\beta_1\in (0,\infty)$ (class of $\FF_p$, $1<p<\infty$)
or $\beta_1=\infty$ (class of $\FF_{\infty}$). Moreover (see \ref{P: l2 inv ME}* below), if $\Gamma$
is ME to a free group and $\beta_1(\Gamma)=p-1$ is an integer,
then $\Gamma\OE\FF_p$.\\
Also, $\FF_{p_1}\!\times  \FF_{p_2}\! \times\! \cdots \times  \FF_{p_j}\ME \FF_{q_1}\!\times \FF_{q_2}\!\times \cdots \times \FF_{q_k}$ i{f}{f} $j=k$ (with $\infty >p_i,q_i\geq 2$)\footnote{The only non zero $\ell^2$ Betti number is the $j$-th, resp. the $k$-th.}, 
and in this case they are commensurable.\\
In the following three claims, the only non zero $\ell^2$ Betti number is in the middle dimension of the associated symmetric space:\\
If two lattices of $\mathrm{Sp}(n,1)$ and $\mathrm{Sp}(p,1)$  are ME, then $p=n$ (see also \ref{P: C(Lambda) is ME inv}).\\
If two lattices of $\mathrm{SU}(n,1)$ and $\mathrm{SU}(p,1)$  are ME, then $p=n$.\\
If two lattices of $\mathrm{SO}(2n, 1)$ and $\mathrm{SO}(2p, 1)$  are ME, then $p=n$.\\
If it would happen that lattices in $\mathrm{SU}(n,1)$ and $\mathrm{SO}(p, 1)$ are ME, then $p=2n$, but we suspect that lattices of $\mathrm{SU}(n,1)$ and $\mathrm{SO}(2n, 1)$ are never ME.

Thanks to the great flexibility in constructing groups with prescribed $\ell^2$ Betti numbers (see \cite{CG86}), it follows from \ref{P: l2 inv ME} that there are uncountably many different ME classes.
\begin{propME}
\item[\theMEpropcourant*]
	If one of the $\beta_n(\Gamma)$ is different from $0,\infty$,
	and if $\Gamma\ME\Lambda$, with index $\iota=[\Gamma:\Lambda]$, 
	then $\iota= \frac{1}{\lambda}$ is imposed. In this case, $\Gamma\OE\Lambda$ 
	i{f}{f} $\lambda=1$ (see \ref{P: ME <-> SOE}).
\end{propME}
\begin{propME}
\item [\theMEpropcourant**]
An infinite group with $\beta_1=0$ is  ME to a free group if and only if it is amenable.
\end{propME}
An amalgamated free product $\Zmath\amalgam{\Zmath}\Zmath$ is ME to a free group
if and only if it is amenable, i.e. the common $\Zmath$ coincides with one the 
components or injects as a subgroup of index 2 in both.

\begin{remark}
If $\Gamma=\Gamma_1\ast_{D}\Gamma_2$ is ME to a free group and nonamenable,
with $D$ infinite, then $\beta_1(\Gamma_1)+\beta_1(\Gamma_2)-\beta_1(D)>0$.
To see it, first observe that $\beta_1(\Gamma)>0$ and $\beta_2(\Gamma)=0$.
The Mayer-Vietoris' exact sequence then gives the estimate.
For instance, amalgamated free products $\FF_2\ast_{\FF_3} \FF_2$ are not ME to $\FF_2$.
\end{remark}

\begin{propME}
\item Ergodic dimension and approximate ergodic dimension are ME invariants 
(see \cite[Sect. 5.3 and Prop. 6.5]{Gab02}).\label{P: ergodic dimension ME inv.}
\end{propME}
Both invariants are useful when all the $\ell^2$ Betti numbers are $0$.
For instance, very little is known about the ME classification 
of lattices in $\mathrm{SO}(m,1)$ for $m$ odd.
However, a reasoning on ergodic dimension shows that:
If two lattices in  $\mathrm{SO}(2n+1,1)$ and $\mathrm{SO}(2p+1,1)$ are ME, with
$p\leq n$, then $n\leq 2p$ \cite[Cor.~6.9]{Gab02}. Similarly, such a group as 
the following $\Gamma$ has approximate dimension $j$, and thus $\Gamma=\FF_{p_1}\!\times  \FF_{p_2}\! \times\! \cdots \times  \FF_{p_j}\times \Zmath\ME\FF_{q_1}\!\times \FF_{q_2}\!\times \cdots \times \FF_{q_k}\times \Zmath $ i{f}{f} $j=k$
(all the $p_i,q_i\geq 2$). In this case they are commensurable. 
\begin{propME}
\item The set $I_{\mathrm{ME}}(\Gamma)$ of all indices of
ME couplings of $\Gamma$
with itself \cite[Sect. 2.2]{Gab00b} is a ME invariant.
\end{propME}
Since non-ergodic ME couplings are allowed, $I_{\mathrm{ME}}(\Gamma)$
is a convex subset of $\Rmath^{*}_{+}$.
If one of the $\ell^2$ Betti numbers $\beta_n(\Gamma)\not=0,\infty$, then
$I_{\mathrm{ME}}(\Gamma)=\{1\}$.
When $\Gamma$ is the direct product of an 
infinite amenable group with any group
then $I_{\mathrm{ME}}(\Gamma)=\Rmath^{*}_{+}$.
\begin{question}
Are there groups $\Gamma$ such that the set of 
all indices of \textbf{ergodic}
ME couplings of $\Gamma$
with itself is discrete $\not=\{1\}$ ?
\end{question}
Recall that the \textbf{fundamental group} of a p.m.p. ergodic $\Gamma$-action $\alpha$ on $(X,\mu)$ is the subgroup of $\Rmath_{+}^{*}$ generated by the set of measures $\mu(A)$ for those $A\subset X$ such that the restriction of $\RR_{\alpha}$ to $A$ is isomorphic with $\RR_{\alpha}$ (equivalently the set of compression constants of $\alpha$ with itself). To compare with, $I_{\mathrm{ME}}$ is related to the bigger set of measures $\mu(A)$ for those $A\subset X$ such that the restriction of $\RR_{\alpha}$ to $A$ may be produced by \emph{some} free action of $\Gamma$.


\medskip
In \cite{MS02}, Monod-Shalom introduce the class $\mathcal{C}_{\mathrm{reg}}$
of all groups $\Gamma$ such that $\mathrm{H}_{b}^{2}(\Gamma,\ell^{2}(\Gamma))$ 
is non-zero and the larger class 
$\mathcal{C}$ 
of all groups $\Gamma$ admitting a
mixing unitary representation $\pi$ such that $\mathrm{H}_{b}^{2}(\Gamma,\pi)$
is non-zero. 
Non amenable free products and non-elementary subgroups of hyperbolic groups 
all belong to the class $\mathcal{C}_{\mathrm{reg}}\subset \mathcal{C}$
 \cite{MS02,MMS03}. A group in the class $\mathcal{C}$ has finite center and 
 is not a direct product of two infinite groups. Also, being in the class 
$\mathcal{C}$ passes to normal subgroups.

\begin{propME}
\item Belonging to the class $\mathcal{C}_{\mathrm{reg}}$ (resp. $\mathcal{C}$ )
is a ME invariant 
\cite{MS02}.
\label{P: C is ME inv}
\end{propME}
For instance, a normal subgroup $\Lambda$ of a group $\Gamma$ ME to a non elementary hyperbolic group has finite center and is not a direct product 
of two infinite groups.
This statement for $\Lambda=\Gamma$ were obtained by S.~Adams 
\cite{Ada94,Ada95}.
Lattices in $\mathrm{SO}(2n,1)$ are not ME to lattices in 
$\mathrm{SL}(2,\Rmath)\times \mathrm{SL}(2,\Rmath)\times\cdots\times \mathrm{SL}(2,\Rmath)$ (follows from \cite{Ada94}).
Observe that hyperbolicity itself is not preserved since $\Zmath^2\ast\Zmath^2\OE\FF_2$.

\subsection{Some ``Rigidity'' Results}
The following result is a spectacular achievement of the rigidity phenomena
attached to higher rank lattices, after Margulis' super-rigidity and Zimmer's 
cocycle super-rigidity \cite{Zim84}. In fact, Zimmer obtained the similar 
result with the additional assumption that the mysterious group admits
a linear representation with infinite image \cite{Zim84,Zim91}.
\begin{propME}
\item (ME Rigidity -- Higher Rank Lattices) Any countable group which is ME to 
	a lattice in a connected simple Lie group $G$ with finite center and real 		rank $\geq 2$, is commensurable up to finite kernel with a lattice in $G$ 		\cite{Fur99a}.
	\label{P: ME rig. latt. high rank}
\end{propME}
The prototype of such Lie groups are $G=\mathrm{SL}(n,\Rmath)$ or
$\mathrm{SO}(p,q)$, for $n-1, p,q\geq 2$.

\medskip

Direct products sometimes appear as analogue to higher rank lattices.
Here, bounded cohomology and belonging to the class $\mathcal{C}$ 
(see \ref{P: C is ME inv}) rigidifies the situation:
\begin{propME}
\item (ME Rigidity -- Products) If 
$\Gamma_1\times\cdots\times\Gamma_n\ME\Lambda_1\times\cdots\times\Lambda_{p}$ 
are ME products of (non-trivial) 
torsion-free countable groups with $n\leq p$, where all the $\Gamma_i$'s are in the class $\mathcal{C}$,
then $n=p$ and after
permutation of the indices 
$\Gamma_i\ME\Lambda_i$ for all $i$ \cite[Th. 1.14]{MS02}.
\label{P: prod gps in C}
\end{propME}
For direct products of free groups, with $p_i,q_i\geq 2$ the value $p_i,q_i=\infty$ being now allowed, 
$\FF_{p_1}\!\times  \FF_{p_2}\! \times\! \cdots \times  \FF_{p_j}\ME
\FF_{q_1}\!\times \FF_{q_2}\!\times \cdots \times \FF_{q_k}$
i{f}{f} $j=k$ and the number of times $\infty$ occurs is the same; 
and in this case they are commensurable (compare with~\ref{P: l2 inv ME}
which gives only $j=k$).
\begin{propME}
\item (ME Rigidity -- Quotients by Radicals) 
Let $M\triangleleft \Gamma$, $N\triangleleft \Lambda$
be amenable normal subgroups of $\Gamma, \Lambda$ such that $\Gamma/M$ and 
$\Lambda/N$ are in the class $\mathcal{C}$ and torsion-free.
If $\Gamma\ME\Lambda$, then $\Gamma/M\ME\Lambda/N$ \cite[Th. 1.15]{MS02}.
\end{propME}

\subsection{Application}

Consider Lie groups of the form
$\prod_{i\in I}\mathrm{Sp}(m_i,1)\times \prod_{j\in J} \mathrm{SU}(n_j,1) \times\prod_{k\in K}\mathrm{SO}(p_k,1)$,
where $I,J,K$ are finite sets and $m_i,n_j,p_k\geq 2$.
Let $\Gamma$ and $\Gamma'$ be lattices of such Lie groups $G$ and $G'$.
Assume they are ME. Each of them is ME to a product of torsion free
cocompact (thus hyperbolic) lattices of the factors $\prod_{t\in I,J,K}\Gamma_t\ME\Gamma\ME\Gamma'\ME\prod_{t'\in I',J',K'} \Gamma_{t'}$.
By \ref{P: prod gps in C}, the pieces correspond under ME, after reordering.
It follows that the number of factors coincide. 
Thanks to property T (\ref{P: T is ME inv}), the pieces in the
$\mathrm{Sp}(\cdot,1)$'s are ME. Moreover, \ref{P: C(Lambda) is ME inv}
or the examination of the $\ell^2$ Betti numbers
(\ref{P: l2 inv ME}) ensures that the sets $\{m_i:i\in I\}$ and $\{m_{i'}:i'\in I'\}$
are the same. Similarly, the number of odd $p_k$, $p_{k'}$ are the same
and the set of $n_j$ and even $p_k$ on the one hand, coincide 
with the set $n_{j'}$ and even $p_{k'}$ on the other hand.
For instance if $K,K'$ are empty, then $\Gamma\ME\Gamma'$
implies that $G,G'$ are isomorphic. And similarly
if $J,J'$ are empty and all the $p_k$ are even.

\section{A Construction of Measure Equivalent Groups}

\subsection{Measure Free-Factor}

We first introduce a measure theoretic notion analogue to
free factors in group theory.
\begin{definition}\label{def: measure free-factor}
A subgroup $\LL<\GG$ is called a \textbf{measure free-factor of $\GG$}
if $\GG$ admits a free p.m.p. action $\alpha$ such that $\RR_{\alpha(\LL)}$
is freely supplemented in $\RR_{\alpha(\GG)}$; i.e. there exists a subrelation $\SS<\RR$
such that $\RR=\RR_{\alpha(\LL)}\ast \SS$. 
\end{definition}
This notion is clearly invariant under the automorphism group of $\GG$.
\\
Example: If $\LL$ is a free-factor of $\GG$, i.e. if $\GG$ decomposes as a free 
product $\GG=\LL\ast \LL'$, then $\LL$ is a measure free-factor.
This definition, motivated by Theorem~\ref{th: iterated amalg fp is ME to free group}, of course challenges to exhibit non-trivial examples.
This is done in the following:
\begin{theorem}\label{th: prod commutators is a measure-free-factor}
The cyclic subgroup $\CC$
generated by the product of commutators $\smash{\kappa:=\prod 
\limits_{i=1}^{i=p}[a_{i},b_{i}]}$ is a measure free-factor in 
the free group $\FF_{2p}=\langle a_{1},\cdots,a_{p},b_{1},\cdots,b_{p}\rangle$.
\end{theorem}
Since it vanishes in the abelianization, $\CC$ is not a free factor in the usual sense.

\medskip
\noindent
\textsc{Proof} of Th.~\ref{th: prod commutators is a measure-free-factor}:
This result is obtained by exhibiting a free p.m.p. action 
$\sigma$ of $\FF_{2p}$ on an $(X,\mu)$ whose equivalence relation is also 
generated by a \textbf{treeing} $\Phi=(\varphi_{\gamma})_{\gamma\in \{\kappa, a_{1},\cdots,a_{p},b_{1},\cdots,b_{p}\} }$ made of an automorphism $\varphi_{\kappa}=\sigma(\kappa)$, defined on the whole of $X$, and of partially 
defined isomorphims $\varphi_{\gamma}$ that are restrictions of 
the generators $\gamma$.

It follows from the classification of surfaces that the free group $\FF_{2p}$ 
is isomorphic with the fundamental group of the oriented surface of genus 
$p$ with one boundary component (supporting the base point), 
via an isomorphism sending the product of commutators $\kappa$ to the boundary curve and the generators to simple
closed curves that are disjoint up to the base point.

The Cayley graph $\GGG$ and the universal cover of the Cayley complex
associated to the (non-free) presentation $\langle 
\kappa, a_{1},\cdots,a_{p},b_{1},\cdots,b_{p}\ \vert\  \kappa:=\prod 
\limits_{i=1}^{i=p}[a_{i},b_{i}]\rangle$
are thus planar, with one orbit of bounded $2$-cells, and boundary edges labelled $\kappa$.
Its dual graph $\GGG^{*}$ 
is a regular tree of valency $2p$
which does not cross the edges labelled $\kappa$.
It is equipped with a natural vertex-transitive and free action of $\FF_{2p}$
which thus may be seen as a ``standard'' Cayley tree for the free group 
$\FF_{2p}$.
Denote by $E^*$ its edge set
and by ${\mathcal{F}^{*}}$ the  subset of $\{0,1\}^{E^*}$, corresponding to the 
(characteristic functions of the) forests, whose connected components are infinite trees 
with one end. The subset ${\mathcal{F}^{*}}$ is equipped with a natural $\FF_{2p}$-action.
\begin{proposition}\label{prop: existence of a one-ended forest}
The subset ${\mathcal{F}^{*}}$ supports an $\FF_{2p}$-invariant probability
measure $\nu$.
\end{proposition}
This is an immediate consequence 
of \cite[Cor~2.6 or Th.~4.2]{Hag98}. For completeness, we will give an elementary proof below; now we continue with the proof of Th.~\ref{th: prod commutators is a measure-free-factor}.

Denote the edge set of $\GGG$ by $E$.
A forest $f^{*}\in {\mathcal{F}^{*}}$ defines a subgraph $H(f^{*})\in 
\{0,1\}^{E}$ by removing the edges of $E$ that $f^{*}$ crosses. 
The end-points of an edge removed (thanks to an edge $e^*\in f^{*}$)
 are connected in $H(f^{*})$ by the finite planar path surrounding 
the finite bush of $f^{*}\setminus\{e^*\}$.
It follows that $H(f^{*})$ is a connected subgraph of $\GGG$ containing 
all the edge labelled $\kappa$ (since the dual $\GGG^*$ does not cross 
the edges of $\GGG$ labelled $\kappa$). By planarity and since all the 
connected components of $f^{*}$ are infinite, 
$H(f^{*})$ has no cycle: it is a tree.

The equivariant map $H:{\mathcal{F}^{*}}\to \{0,1\}^{E}$ thus pushes the measure $\nu$ to an $\FF_{2p}$-invariant probability measure $\mu$, supported on the set of (connected) subtrees of $\GGG$, containing all the edges labelled $\kappa$.

Making the action free if necessary (by considering the diagonal action on $(X,\mu)=(X',\mu')\times (\{0,1\}^{E},\mu)$ for some free p.m.p. $\FF_{2p}$-space $(X',\mu')$, see \cite[1.3.c]{Gab04}), and according, for example, to \cite[1.3.f]{Gab04} this is equivalent with the existence of the required graphing on $(X,\mu)$. More precisely, the natural $\FF_{2p}$-equivariant map $\pi: X\to \{0,1\}^{E}$ forgetting the first coordinate sends $x\in X$ to the subset $\pi(x)\subset E$ of the edge set of $\GGG$. For $\gamma\in \{\kappa, a_{1},\cdots,a_{p},b_{1},\cdots,b_{p}\}$ define the Borel subset $A_{\gamma}$ of those $x\in X$ for which the edge $[id,\gamma]$ of $\GGG$ belongs to $\pi(x)$ and consider the partially defined isomorphism $\varphi_{\gamma}:=\gamma^{-1}_{\vert A_{\gamma}}$. The graphing $\Phi=(\varphi_{\kappa},\varphi_{a_{1}},\cdots,\varphi_{a_{p}},\varphi_{b_{1}},\cdots,\varphi_{b_{p}})$ matches the conditions:
The graph associated with the orbit of $x$ is isomorphic with the subgraph defined by $\pi(x)$, and thus it is a tree for $\mu$-almost all $x\in X$ ($\Phi$ is a treeing), $\Phi$ generates the same equivalence relation as the $\FF_2$-action since the graph $\pi(x)$ is connected, and $A_{\kappa}=X$ (up to a $\mu$-null set), i.e. $\varphi_{\kappa}$ is defined almost everywhere.
\endofproof

\bigskip\noindent\textsc{Proof} of Prop.~\ref{prop: existence of a one-ended forest}:
Let $(Y,m)$ be a Borel standard probability space with a partition $Y=\coprod_{i\in \Nmath} (U_i \coprod V_i)$, where $m(U_i)=m(V_i)=\frac{1}{2^{i+1}}$. Let $\hhh_{1},\cdots,\hhh_{2p}$ be a family of p.m.p. automorphisms of $(Y,m)$ such that $\hhh_1(U_i)=U_{i+1}\cup V_{i+1}$ and $\hhh_2(V_i)=U_{i+1}\cup V_{i+1}$. They define a p.m.p. (non-free) action $\alpha$ of $\FF_{2p}=\langle \gamma_{1},\cdots,\gamma_{2p}\rangle$ on $Y$. Denote $A_1=\coprod_{i\in \Nmath} U_i$ and $A_2=\coprod_{i\in \Nmath} V_i$. The restrictions of the automorphisms $\varphi_{1}=\hhh_{1\vert A_1}$ and $\varphi_{2}=\hhh_{2\vert A_2}$ define a graphing, which is a treeing whose associated graphs are all one-ended trees, and which generates a subrelation of $\RR_{\alpha}$. If $e_{i}^{*}$ denotes the edge $[id,\gamma_i^{-1} id]$ in the Cayley graph of $\FF_{2p}$ associated with the generating family $(\gamma_i)$, then the conditions $\pi(y)(e_{1}^{*})=1$ i{f}{f} $y\in A_1$,  $\pi(y)(e_{2}^{*})=1$ i{f}{f} $y\in A_2$, and $\pi(y)(e_{i}^{*})=0$ for the other $i$'s, extend by $\FF_{2p}$-equivariance to a map $Y\to {\mathcal{F}^{*}}$. Pushing forward the invariant measure $m$ delivers an instance of the required measure $\nu$.\endofproof
\begin{proposition}
If $\Lambda<\Gamma$ is a measure free-factor of $\Gamma$ and $\Gamma'$ is any countable group, then $\Lambda$ as well as $\Lambda\ast \Gamma'$ is a measure free-factor of $\Gamma\ast\Gamma'$.
\end{proposition}
\noindent\textsc{Proof}: Consider an action of $\Gamma$, that witnesses the measure free-factor condition. Look at it as an action of $\Gamma\ast \Gamma'$ via the obvious map $\Gamma\ast \Gamma'\to \Gamma$, and consider the direct product action $\alpha$ with any free p.m.p. action of $\Gamma\ast \Gamma'$. The lifting properties of locally one-to-one morphisms (see Lemma~\ref{lem: loc 1-to-1 and splitting}) applied to the restriction of $\alpha$ to $\Gamma$ show that $\alpha$ matches the required condition.\endofproof

\begin{corollary}\label{cor: F4 measure free-factor of Fq}
Let $\Sigma$ be an oriented surface with $r>0$ boundary components. Let $\pi_1(\Sigma,*)\simeq \FF_{q}$ be its fundamental group and $\Lambda\simeq \FF_{r}$ the subgroup generated by the boundary components. Then $\Lambda$ is a measure free-factor of $\FF_q$.
\end{corollary} 
\textsc{Proof:} $\pi_1(\Sigma,*)$ admits a free generating set $a_i,b_i,\kappa_1,\kappa_2,\cdots \kappa_{r-1}$ for which $\prod_{i=1}^{p}[a_i,b_i], \kappa_1,\kappa_2,\cdots \kappa_{r-1}$ freely generates $\Lambda$. In fact, the product of the boundary components $\kappa_1,\kappa_2,\cdots \kappa_{r}$ equals $\prod_{i=1}^{p}[a_i,b_i]$.\endofproof

\subsection{Application}

\begin{theorem}\label{th: iterated amalg fp is ME to free group}  If $\GG$ is ME to the free group $\FF_p$ ($p=2$ or $p=\infty$) and admits a measure free-factor subgroup $\LL$, then the iterated amalgamated free product $\itFP{\LL}{n}\GG=\GG\amalgam{\LL}\GG\amalgam{\LL}\cdots \amalgam{\LL}\GG$ is ME to $\FF_p$,  and the infinite iterated amalgamated free product ${\itFP{\LL}{\infty}\GG=\GG\amalgam{\LL}\GG\amalgam{\LL}\cdots \amalgam{\LL}\GG\amalgam{\LL}\cdots}$ is ME~to~$\FF_{\infty}$. Moreover, $\LL$ remains a measure free-factor in the resulting group.
\end{theorem}

Applied with Th.~\ref{th: prod commutators is a measure-free-factor}, this gives Th.~\ref{th: branched surfaces ME F2} of the introduction:
\begin{corollary}\label{cor: branched surfaces ME F2} For each $n\in \Nmath$, the iterated amalgamated free product $\itFP{\CC}{n}\FF_{2p}$ is measure equivalent to the free group $\FF_2$. The infinite amalgamated free product: $\itFP{\CC}{\infty}\FF_{2p}$ is measure equivalent to the free group $\FF_{\infty}$. Moreover,  $\CC$ is a measure free-factor in the resulting group.
\end{corollary}

\textsc{Proof} of Th.~\ref{th: iterated amalg fp is ME to free group}: The diagonal action $\beta$ on the product measure space $X$ of the two actions witnessing each condition on $\GG$, besides being free and p.m.p.,both is treeable (Lemma~\ref{lem: loc 1-to-1 -> graphing}) and admits a splitting as $\SS=\SS_{\beta(\LL)}\ast \SS'$ (the splitting condition lifts, see Lemma~\ref{lem: loc 1-to-1 and splitting}). It thus admits a treeing  that splits into two pieces $\Phi=\Phi_{\LL}\vee\Phi'$ where $\Phi_{\LL}$ generates $\SS_{\beta(\LL)}$ and $\Phi'$ generates $\SS'$ ($\SS'$ is treeable by \cite[Th. 5]{Gab00}). Consider\\
(1) the natural surjective homomorphism $\pi:\itFP{\LL}{n}\GG\to \GG$, defined by the identity on each copy~of~$\GG$,\\
(2) the (non free) p.m.p. action of $\itFP{\LL}{n}\GG$ on $X$ via $\beta\circ\pi$,\\
(3) any free p.m.p. action of $\itFP{\LL}{n}\GG$ on a standard Borel space $Y$,\\ 
(4) the diagonal $\itFP{\LL}{n}\GG$-action $\alpha^{n}$ on the product measure space $Z:=X\times Y$ (it is p.m.p.~and~free),\\
(5) the $\itFP{\LL}{n}\GG$-equivariant projection to the first coordinate $\Pi:X{\times}Y{\to}X$.

\medskip
For each $i=1,\cdots,n$, denote by $\RR_i:=\RR_{\alpha^{n}(\GG_i)}$ the equivalence relation generated by the restriction of the action $\alpha^{n}$ to the $i$-th copy $\GG_i$ of $\GG$ (resp.  $\RR_{\LL}:=\RR_{\alpha^{n}(\LL)}$ for the restriction to $\LL$). The morphism $\Pi$ is locally one-to-one when restricted to $\RR_i$ (Ex.~\ref{ex: p loc 1-to-1}). Denote by $\tilde{\Phi}_{i}=\tilde{\Phi}_{i\;\LL}\vee\tilde{\Phi}'_{i}$ the graphing of $\RR_i$ given from $\Phi$ by the lifting lemma~\ref{lem: loc 1-to-1 -> graphing}. 
More concretely: Up to subdividing the domains, each generator $\varphi:A_{\varphi}\to B_{\varphi}$ of the treeing $\Phi$ is the restriction of the $\beta$-action of one element 
$\gamma_\varphi$ of $\GG$. Each $\sigma_j$ of the $n$ homomorphic sections $\sigma_j:\GG\to \itFP{\LL}{n}\GG$ of $\pi$ delivers a graphing on $Z$ of the following form: $\tilde{\Phi}_{j}=\tilde{\Phi}_{j\LL}\vee\tilde{\Phi}'_{j}$, where each element $\tilde{\varphi}:\tilde{A_{\varphi}}\to\tilde{B_{\varphi}}$ is defined from the element $\varphi A_{\varphi}\to B_{\varphi}$ of $\Phi$ as the restriction of the $\alpha^{n}$-action of $\sigma_j(\gamma_{\varphi})$ to $\tilde{A_{\varphi}}:=\Pi^{-1}({A_{\varphi}})$.

By the coincidence of the copies of $\LL$, it follows that all the $\tilde{\Phi}_{i\;\LL}$ coincide, so that $\tilde\Phi=\tilde{\Phi}_{1\;\LL}\vee \tilde{\Phi}'_{1}\vee\tilde{\Phi}'_{2}\cdots\vee \tilde{\Phi}'_{n}$ still generates $\RR_{\alpha^{n}}$. Since it lifts a treeing, each $\tilde{\Phi}'_{i}$ (as well as $\tilde{\Phi}_{1\;\LL}$) is itself a treeing. The structure of free product with amalgamation of ${\RR=\RR_{1} \amalgam{\RR_{\LL}}\RR_{2}\amalgam{\RR_{\LL}} \cdots \amalgam{\RR_{\LL}}\RR_{n}}$ ensures that $\tilde\Phi$ is globally a treeing: $\itFP{\LL}{n}\GG$ admits a treeable free action, and $\RR_{\LL}$ is a free-factor. A similar reasoning give the similar result for $\itFP{\LL}{\infty}\GG$. These groups are consequently ME to a free group by \ref{P: Hjorth result}. The precise ME-class is determined by the first $\ell^2$ Betti number (by \ref{P: l2 inv ME}).
\endofproof

\begin{remark}
The treeing $\Phi$ produced in the proof of 
theorem~\ref{th: prod commutators is a measure-free-factor} 
cannot be a quasi-isometry, 
i.e. there is some generator $\gamma\in \{a_{1},\cdots,a_{p},b_{1},\cdots,b_{p}\}$
for which the length of the $\Phi$-words sending $x$ to $\gamma(x)$
is not bounded a.s. For otherwise, the graphing produced in the 
proof of theorem~\ref{th: iterated amalg fp is ME to free group}
for $n=2$ would equip each equivalence class with a graph structure (tree)  quasi-isometric to the group itself $\FF_{2p}\ast_{\CC}\FF_{2p}$ (the fundamental group of a closed compact surface, without boundary). 
It follows from the particular form of $\Phi$ that the surface group 
$\langle a_1, \cdots, a_p, b_1, \cdots, b_p\rangle \ast_{C=C'}
\langle a'_1, \cdots, a'_p, b'_1,\cdots, b'_p\rangle$
has a treeing made of $\CC$ and restrictions of the ``natural'' generators
$a_i, b_i, a'_i, b'_i$.
\end{remark}

\begin{question}
What are all the measure free-factors of the free group $\FF_2$~?
\end{question}
Observe that in a group $\Gamma$ ME to a free group, 
a cyclic subgroup $D$ strictly contained in a greater cyclic subgroup $E$
cannot be a measure free-factor by Th.~\ref{th: iterated amalg fp is ME to free group}, since iterated amalgamated free product
$\Gamma\amalgam{D} \Gamma\amalgam{D} \Gamma$  is not ME to a free group: 
it contains the nonamenable subgroup $E\amalgam{D} E \amalgam{D} E$ for which $\beta_1=0$, then use \ref{P: l2 inv ME}** and \ref{P: ME F stable subgroup}.
\begin{question}
It is known that if an amalgamated free product 
$\FF_p\ast_{\Zmath}\FF_q$ happens to be a free group, then $\Zmath$
is a free factor in one of $\FF_p$ or $\FF_q$ (see \cite[Ex. 4.2]{BF94}). 
Is it true that similarly if $\FF_p\ast_{\Zmath}\FF_q\ME\FF_2$ then $\Zmath$ is a measure free-factor in $\FF_p$ or $\FF_q$ ?
\end{question}

\subsection{Orbit Equivalence of Pairs}
We will now make use of the following refined notion of Orbit Equivalence.
\begin{definition}
Consider for $i=1,2$ a group $\Gamma_i$ and a subgroup $\Gamma_i^0$.
An \textbf{orbit equivalence of the pairs} $(\Gamma_1^0<\Gamma_1)$
and $(\Gamma_2^0<\Gamma_2)$ is the data of two p.m.p. actions
$\alpha_1$ and $\alpha_2$ of $\Gamma_1$ and $\Gamma_2$
on the probability measure standard Borel space $(X,\mu)$, that generate the same equivalence relation
($\RR_{\alpha_1(\Gamma_1)}=\RR_{\alpha_2(\Gamma_2)}$)
and such that the restrictions of $\alpha_{1}$ to $\Gamma_1^0$
and of $\alpha_{2}$ to $\Gamma_2^0$ also define a common sub-relation:
($\RR_{\alpha_1(\Gamma_1^0)}=\RR_{\alpha_2(\Gamma_2^0)}$).
The existence of such an equivalence of pairs is denoted by:
$$(\Gamma_1^0<\Gamma_1)\OE(\Gamma_2^0<\Gamma_2)$$
\end{definition}
\begin{definition} It is a \textbf{strong orbit equivalence of the pairs} 
if moreover 
the sub-groups $\Gamma_1^0$ and $\Gamma_2^0$ are isomorphic, via
an isomorphism $\phi: \Gamma_1^0\to \Gamma_2^0$ that turns the actions 
conjugate:
$\forall \gamma\in \Gamma_1^0, \alpha_1(\gamma)=\alpha_2(\phi(\gamma))$
in $\mathrm{Aut}(X,\mu)$. The existence of such an equivalence of pairs is denoted by:
$$(\Gamma_1^0<\Gamma_1)\OEstrong(\Gamma_2^0<\Gamma_2)$$
\end{definition}
Measure free-factors naturally lead to orbit equivalence of pairs:
\begin{theorem}\label{th: examples of OE pairs}
Let $\GG$ be a group ME to the free group $\FF_2$ and let $\GG_0$
be a subgroup such that $\beta_1(\GG)-\beta_1(\GG_0)=q$ is an integer.
If $\GG_0$ is a measure free-factor of 
$\GG$, then there exists a strong orbit equivalence of pairs
$(\GG_0<\GG)\OEstrong(\GG_0<\GG_0\ast \FF_{q})$.
\end{theorem}
\noindent
\textsc{Proof} of Th. \ref{th: examples of OE pairs}:
Like in the proof of Th.~\ref{th: iterated amalg fp is ME to free group},
consider (by taking a diagonal action of two actions witnessing both
properties of $\GG$)
a free p.m.p. action $\beta$ of $\GG$ that  both is treeable and
admits a splitting $\SS_{\beta}=\SS_{\beta(\GG_0)}\ast \SS'$.
Up to considering the ergodic decomposition, $\beta$ may be assumed to be ergodic.
The subrelation $\SS_{\beta(\GG_0)}$ admits a treeing $\Phi$
and  $\SS'$ admits a treeing $\Psi$ of cost $q$ 
\cite[Cor. 3.16, Cor. 3.23]{Gab02}). A recursive use of Theorem 28.3 of \cite{KM04} (due to Hjorth) gives that $\SS'$ may be replaced in the above decomposition by $\SS''$ produced by a 
free action of the free group $\FF_{q}$. The free product decomposition of $\SS_{\beta}=\SS_{\beta(\GG_0)}\ast \SS''$ 
asserts that the action of $\GG_0$ fits well to produce the pair $\SS_{\beta(\GG_0)}<\SS_{\beta}$
by a free action of $\GG_0\ast \FF_{q}$.
The resulting orbit equivalence of pairs is strong by construction.
\endofproof 
\begin{example}\label{ex: examples of OE pairs 1}
If $\CC$ is the cyclic subgroup generated by the product of commutators $\kappa$ in $\FF_{2p}$
then $(\CC<\FF_{2p})\OEstrong(\Zmath<\Zmath\ast \FF_{2p-1})$.
\end{example}

\begin{example}\label{ex: examples of OE pairs 2}
If $\Lambda$ is the subgroup of $\FF_{q}$ as in 
Corollary~\ref{cor: F4 measure free-factor of Fq}, then 
$(\FF_r\simeq\Lambda<\FF_q)\OEstrong(\FF_r<\FF_r\ast \FF_{q-r})$.
\end{example}

\begin{remark}\label{rem: free prod -> strong}
Observe that  in a free product situation
$(\Gamma^0<\Gamma_1)\OE(\Gamma^0<\Gamma^0\ast \Gamma'_2)$,
the orbit equivalence of pairs 
may always be assumed strong, by changing the $\Gamma^0$-action
in the free product.
\end{remark}

\begin{theorem}\label{th: amalg+coind}
If $(\GG_1^0<\GG_1)\smash{\OEstrong}(\GG_2^0<\GG_2)$ are two pairs of groups 
admitting a strong orbit equivalence of pairs
and  $G$ is a countable group with a subgroup $\GG^0$ isomorphic with
$\GG_1^0$ and $\GG_2^0$, then
the following pairs 
admit a strong orbit equivalence of pairs:
$$(G<G\amalgam{\GG^0=\GG_1^0}\GG_1)\OEstrong(G<G\amalgam{\GG^0=\GG_2^0}\GG_2).$$
In particular, if $(\Gamma^0<\Gamma_1)\OE(\Gamma^0<\Gamma^0\ast \Gamma'_2)$
and $G$ contains a subgroup isomorphic with $\Gamma^0$, then
$$(G<G\amalgam{\Gamma^0} \Gamma_1) \OEstrong (G<G\ast \Gamma'_2).$$
\end{theorem}
\begin{corollary}\label{cor: f.p. amalg. ME free group}
Let $G$ be any countable group and $H$  an infinite cyclic subgroup.
Let $\CC$ be the cyclic subgroup generated by the product of commutators 
$\kappa$ in $\FF_{2p}$. Then 
$\smash{G\amalgam{H=\CC}\FF_{2p}\OE G\amalgam{}\FF_{2p-1}}$.
In particular, if $G$ is ME to $\FF_2$, then $G\amalgam{H=\CC}\FF_{2p}\ME\FF_2$.
\end{corollary}

\noindent
\textsc{Proof} of Th.~\ref{th: amalg+coind}:
Let $\beta_1$ and $\beta_2$ be p.m.p. actions on $(Y,\nu)$ defining a strong
OE of the pairs given by the assumption. 
Consider a free p.m.p. action $\tilde\alpha_1$ of $G\amalgam{\GG^0=\GG_1^0}\GG_1$ on some $(X,\mu)$, with a $\GG_1$-equivariant 
map $p:Y\to X$ (the co-induced action from $\GG_1$ to $G\amalgam{\GG^0=\GG_1^0}\GG_1$ satisfies this property, see section~\ref{subsect: co-ind action} below).
Call $\alpha_1$ the restriction of $\tilde\alpha_1$ to $\GG_1$.
The locally one-to-one morphism $p$ from $\RR_{\alpha_1(\GG_1)}$
to $\SS_{\beta_1(\GG_1)}=\SS_{\beta_2(\GG_2)}$ allows to lift
(cf. Lemma~\ref{lem: loc 1-to-1 -> action})
the action $\beta_2$ to an action $\alpha_2$ of $\GG_2$ on $Y$, 
which is orbit equivalent with $\alpha_1$. By uniqueness in 
Lemma~\ref{lem: loc 1-to-1 -> action}, strongness also lifts, 
i.e. $\alpha_1$ and $\alpha_2$ coincide on $\GG_1^0$ and
$\GG_2^0$, via the given isomorphism between these groups\footnote{Thus, 
$\alpha_1$ and $\alpha_2$ also form a strong orbit equivalence for the pairs
$(\GG_1^0<\GG_1)$ and $(\GG_2^0<\GG_2)$; with an additional potentiality.}.
Now, using $\tilde\alpha_1$, one extends the $\GG_2^0$-restriction of $\alpha_2$
to a free action of $G$, so as to produce an action $\tilde\alpha_2$ of ${G\amalgam{\GG^0=\GG_2^0}\GG_2}$.
Given the amalgamated free product structure of $\RR_{\tilde\alpha_2}=\RR_{\tilde\alpha_2(G)}\amalgam{\RR_{{\tilde\alpha}_2 (\GG^0=\GG_1^0)}}\RR_{\tilde{\alpha}_2(\GG_1)}$, 
and since $\tilde\alpha_1$ and $\tilde\alpha_2$ produce the same equivalence relation when restricted, on the one hand to $G$, on the second hand to $\GG_1^0$ and $\GG_2^0$, and on the third hand (!) to $\GG_1$ and $\GG_2$, it follows that
the action $\tilde\alpha_2$ is free, and forms with 
$\tilde\alpha_1$ the required strong orbit equivalence of pairs. 
For the ``in particular'' part, just observe 
with remark~\ref{rem: free prod -> strong}
that $(G<G\amalgam{\Gamma^0} \Gamma_1) \OEstrong(G<G\amalgam{\Gamma^0}
\Gamma^0\ast \Gamma'_2)=(G<G\ast \Gamma'_2)$.
\endofproof
\begin{question}
A limit group is a finitely generated group $\Gamma$ that is $\omega$-residually free, i.e. for every finite subset $K\subset \Gamma$ there exists a homomorphism $\Gamma\to F$ to a free group, that is injective on $K$.
It is a natural question to ask whether limit groups are ME to a free group.
\end{question}

\subsection{Co-induced action}
\label{subsect: co-ind action}

Let $\AAA$ be a countable group and $\BBB$ a subgroup. The co-induction is a canonical way, from an action $\beta$ 
of the ``small'' group $\BBB$ on $\YYY$, to produce an action $\alpha$ of the overgroup $\AAA>\BBB$ on a space 
$\XXX$ together with a surjective $\BBB$-equivariant map $\XXX\to \YYY$. Let $\lambda$ and $\rho $ be the left $\AAA$-actions by left multiplication and right multiplication by the inverse on $\AAA$: $\lambda(g) : h\mapsto gh$ and $\rho(g) : h\mapsto hg^{-1}$.

\noindent
The space $\XXX:=\mathrm{coInd}_{\BBB}^{\AAA} \YYY$ is the set of $(\BBB,\rho, \beta)$-equivariant maps from $\AAA$ to $\YYY$: $$\mathrm{coInd}_{\BBB}^{\AAA} \YYY:=\{\Psi:\AAA\to \YYY \ |\ \Psi(\rho(b)a)=\beta(b)(\Psi(a)); \forall b\in \BBB, a\in \AAA\}.$$ The \textbf{co-induced action} $\alpha:=\mathrm{coInd}_{\BBB}^{\AAA} \beta$ is the action of $\AAA$ on $\mathrm{coInd}_{\BBB}^{\AAA} \YYY$ defined from the $\AAA$-action $\lambda$ on the source: $$\mathrm{coInd}_{\BBB}^{\AAA} \beta(\aaa): \Psi(.)\mapsto\Psi(\lambda(\aaa)(.)), \ \ \forall \aaa\in \AAA$$ The map $\mathrm{coInd}_{\BBB}^{\AAA} \YYY\to \YYY ; \Psi\mapsto \Psi(1),$ where $1$ is the identity element of $\AAA$, is $\BBB$-equivariant.

A section $s: \AAA/\BBB\to \AAA$ of the canonical map being chosen, the $\AAA$-space $\mathrm{coInd}_{\BBB}^{\AAA} \YYY$ naturally identifies with $\YYY^{\AAA/\BBB}=\prod_{h\in \AAA/\BBB} \YYY_h$, equipped with an action $\sigma$ of $\AAA$ which permutes the coordinates $\YYY_h$ and then permutes the points in each coordinate $\YYY_h$, via the $\beta$ action of an element of $\BBB$. More precisely, for all $\aaa\in \AAA,\ f\in \YYY^{\AAA/\BBB},\  h\in \AAA/\BBB$, the action gives $[\sigma(\aaa)(f)]_h=\beta(\bbb^{-1})(f_{h'})$, where $\bbb\in\BBB$ and $h'\in \AAA/\BBB$ are defined by the equation $s(h')\bbb=\aaa^{-1} s(h)$.

It follows that if $\beta$ preserves a probability measure $\nu$ on $\YYY$, then $\sigma$ preserves the product probability measure on $\YYY^{\AAA/\BBB}$, and in view of the description of $\sigma$, the corresponding $(\AAA,\alpha)$-invariant measure $\mu$ on $\XXX$ is independent of the choice of the section $s$. If $\beta$ is essentially free, then $\alpha$ is essentially free.

Remark: It is interesting to compare co-induction with the usual notion of induction (the $A$-action induced by $\lambda$ on the second factor, on the quotient of $Y\times A$ by the diagonal $B$-action $(\beta, \rho)$). They produce respectively right-adjoint and left-adjoint functors for the restriction functor $res: \mathcal {A}\to \mathcal{B}$ between the categories $\mathcal{A}$ of spaces with an $A$-action and $\mathcal{B}$ of spaces with a $B$-action, with equivariant maps as morphisms.
In our context above, induction delivers an action of $A$ on $(Y\times A)/(B,\beta, \rho)\simeq Y\times A/B$ with a natural invariant measure which is \emph{infinite} when $B$ has infinite index in $A$.

\bigskip
\nobreak
\noindent \textsc{D.~G.: UMPA, UMR CNRS 5669, ENS-Lyon, 
69364 Lyon
Cedex 7, FRANCE}

\noindent \texttt{gaboriau@umpa.ens-lyon.fr}

\end{document}